\documentclass[12pt]{article}
\usepackage[margin=2.5cm]{geometry}
\newtheorem{Th}{Theorem}
\newtheorem{Sl}{Corollary}
\newtheorem{Statement}{Statement}

\newcommand{\conv}{\mathop{\rm conv}\nolimits}
\newcommand{\Vol}{\mathop{\rm Vol}\nolimits}

\title{\bf On the $1$-convexity of random points in~the~$d$-dimensional spherical layer}
\author{S.\,V.\,Sidorov \\    
{Nizhni Novgorod State University, sesidorov@yandex.ru}}
\date{}

\begin{document}

\maketitle

\abstract{We consider the set of points chosen randomly, independently and uniformly in the $d$-dimensional spherical layer. A set of points is called $1$-convex if all its points are vertices of the convex hull of this set. In \cite{3} an estimate for the cardinality of the set of points for which this set is $1$-convex with probability close to $1$ was obtained.
In this paper we obtain an improved estimate.}

\section{Introduction}

Following \cite{2} we call a finite set of points $\{Y_1,\ldots,Y_k\}\subset\mathbf{R}^d$ $1$-convex if the set of vertices of the convex hull $\conv(Y_1,\ldots,Y_k)$  coincides with $\{Y_1,\ldots,Y_k\}.$ 

Let $B_d=\{x\in\mathbf{R}^d:~ |x|\leq 1\}$ be the $ d $ -dimensional unit ball with center at the origin and $rB_d$ be the $d$-dimensional ball of radius $ r $ with center at the origin. Let $B_d\setminus rB_d$ be the spherical layer and choose points $M_n=\{X_1,\ldots,X_n\}\subset B_d\setminus rB_d$ randomly, independently and according to the uniform distribution on $B_d\setminus rB_d.$ Uniform distribution means that the probability of choosing a random point $ X $ from some $S\subseteq B_d\setminus rB_d$  is proportional to the $d$-dimensional volume of $S$ namely 
$$
P(X\in S)=\frac{\Vol(S)}{\Vol(B_d\setminus rB_d)}.
$$ 
Denote by $A_n$ the event that 
the set $M_n$ is $1$-convex.

In the recent paper \cite{3} it is shown that for all values $0<r<1$, $0<\alpha<1$, 
$$
n<\left(\frac{r}{\sqrt{1-r^2}}\right)^d\left(\sqrt{1+\frac{2\alpha(1-r^2)^{d/2}}{r^{2d}}}-1\right)
\eqno(1)
$$ 
the inequality $P(A_n)>1-\alpha$ holds.

\section{Main result}

The following theorem gives an improved estimate for the number of points $n$ guaranteeing the 1-convexity of a random $n$-element set of points with probability at least $1-\alpha.$ The proof uses an approach used in \cite{2}.

\begin{Th} 
Let $0\leq r<1,$ $0<\alpha<1,$  $n<\sqrt{\alpha 2^{d}(1-r^d)}.$ Then $P(A_n)>1-\alpha.$
\end{Th}
\textsc{Proof.} Denote the event that $X_i\notin \conv(M_n\setminus\{X_i\})$ by $C_i$ ($i=1,\ldots,n$). Clearly $A_n=C_1\cap\ldots\cap C_n$ and $P(A_n)=P(C_1\cap\ldots\cap C_n)=
1-P(\overline{C_1}\cup\ldots\cup \overline{C_n})\geq 1-\sum\limits_{i=1}^{n}P(\overline{C_i}).$ 
Let us find the upper bound for the probability of the event $\overline{C_i}.$ This event means the point $ X_i $ belong to the convex hull of the remaining points, that is $X_i\in \conv(M_n\setminus\{X_i\}).$ Since points in $ M_n $ have the uniform distribution, the probability of $\overline{C_i}$ is equal to 
$$
P(\overline{C_i})=\frac{\Vol(\conv(M_n\setminus\{X_i\}))-\Vol(\conv(M_n\setminus\{X_i\})\cap rB_d)}{\Vol(B_d)-\Vol(rB_d)}
$$ for  $i=1,\ldots,n.$

Hence 
$$
P(\overline{C_i})\leq\frac{\Vol(\conv(M_n\setminus\{X_i\}))}{\gamma_d (1-r^d)},
$$ 
where 
$\Vol(rB_d)=\gamma_d r^d$ is the volume of a ball of radius $r$.

It remains to estimate $\Vol(\conv(M_n\setminus\{X_i\}))$. 
Let 
$$
V(k,d)=\sup \{ \Vol(\conv(Y_1,\ldots,Y_k)):~ Y_1,\ldots,Y_k\in B_d \}
$$ 
and 
$$
V_r(k,d)=\sup \{ \Vol(\conv(Y_1,\ldots,Y_k)):~ Y_1,\ldots,Y_k\in B_d\setminus rB_d \}.
$$ 
Clearly, $V_r(k,d)\leq V(k,d)$ for $0<r<1.$ 
It is known (see e.g. \cite{1}) that 
$$
V(k,d)\leq\frac{k\gamma_d}{2^d}.
$$ 
Since 
$$
\Vol(\conv(M_n\setminus\{X_i\}))\leq V_r(n-1,d)\leq V(n-1,d)\leq\frac{(n-1)\gamma_d}{2^d}
$$ 
then 
$$
P(\overline{C_i})\leq\frac{n-1}{2^d(1-r^d)} \qquad (i=1,\ldots,n).
$$ 
Hence 
$$
P(A_n)\geq 1-\sum\limits_{i=1}^{n}P(\overline{C_i})\geq 1-\frac{n(n-1)}{2^d(1-r^d)}\geq 1-\frac{n^2}{2^d(1-r^d)}.
$$ 
Thus if $n$ satisfies the condition $\frac{n^2}{2^d(1-r^d)}<\alpha$, i.e. $n<\sqrt{\alpha 2^{d}(1-r^d)}$, then the inequality $P(A_n)>1-\alpha$ holds. The theorem is proved. 
\medskip

Let us compare the obtained bound $n<\sqrt{\alpha 2^{d}(1-r^d)}$ with the bound (1)
proposed in \cite{3}.

\begin{Statement} 
Let $g=\left(\frac{r}{\sqrt{1-r^2}}\right)^d\left(\sqrt{1+\frac{2\alpha(1-r^2)^{d/2}}{r^{2d}}}-1\right),$ $0<r<1,$ $0<\alpha<1,$ $d\in\mathbf{N}.$ If $r$ and $\alpha$ are fixed then the following asymptotic estimates hold:

\begin{enumerate}
	\item[\rm 1.] $g\sim\frac{\alpha}{r^d},$ if  $\sqrt{\frac{\sqrt{5}-1}{2}}<r<1.$

	\item[\rm 2.] $g\sim\frac{2\alpha}{r^d(\sqrt{1+2\alpha}+1)}=
\frac{\sqrt{1+2\alpha}-1}{r^d}=(\sqrt{1+2\alpha}-1)(\frac{\sqrt{5}+1}{2})^{d/2},$ if $r=\sqrt{\frac{\sqrt{5}-1}{2}}.$

	\item[\rm 3.] $g\sim\frac{\sqrt{2\alpha}}{(1-r^2)^{d/4}},$ if $0<r<\sqrt{\frac{\sqrt{5}-1}{2}}.$
\end{enumerate}
\end{Statement}

\textsc{Proof.} 
We have $$g=\frac{\left(\frac{r}{\sqrt{1-r^2}}\right)^d\frac{2\alpha(1-r^2)^{d/2}}{r^{2d}}}{\sqrt{1+\frac{2\alpha(1-r^2)^{d/2}}{r^{2d}}}+1}=
\frac{2\alpha}{r^d\left(\sqrt{1+2\alpha\left(\frac{\sqrt{1-r^2}}{r^{2}}\right)^d}+1\right)}.$$ 

If $0<\frac{\sqrt{1-r^2}}{r^{2}}<1$ then $g\sim\frac{\alpha}{r^d}$.

If $\frac{\sqrt{1-r^2}}{r^{2}}=1$ then $g\sim\frac{2\alpha}{r^d(\sqrt{1+2\alpha}+1)}=
\frac{\sqrt{1+2\alpha}-1}{r^d}$.

If $\frac{\sqrt{1-r^2}}{r^{2}}>1$ then $g\sim\frac{2\alpha}{r^d\sqrt{2\alpha\left(\frac{\sqrt{1-r^2}}{r^{2}}\right)^d}}=\frac{\sqrt{2\alpha}}{(1-r^2)^{d/4}}$.

The equality $\frac{\sqrt{1-r^2}}{r^{2}}=1$ holds if $r^4+r^2-1=0,$ that is $r^2=\frac{\sqrt{5}-1}{2},$ $r=\sqrt{\frac{\sqrt{5}-1}{2}}.$ The inequality $0<\frac{\sqrt{1-r^2}}{r^{2}}<1$ holds if  $r^4+r^2-1>0,$ that is for $\sqrt{\frac{\sqrt{5}-1}{2}}<r<1.$ The inequality $\frac{\sqrt{1-r^2}}{r^{2}}>1$ holds if   $r^4+r^2-1<0,$ that is for $0<r<\sqrt{\frac{\sqrt{5}-1}{2}}.$
The statement is proved.

\begin{Sl} 
Let $f=\sqrt{\alpha 2^{d}(1-r^d)},$ $g=\left(\frac{r}{\sqrt{1-r^2}}\right)^d\left(\sqrt{1+\frac{2\alpha(1-r^2)^{d/2}}{r^{2d}}}-1\right),$ $0<r<1,$ $0<\alpha<1,$ $d\in\mathbf{N}.$ If $r$ and  
$\alpha$ are fixed then the following asymptotic estimates of the quotient $\frac{f}{g}$ hold:

\begin{enumerate}
	\item[\rm 1.] $\frac{f}{g}\sim\frac{1}{\sqrt{\alpha}}(r\sqrt{2})^d\rightarrow\infty,$ if  $\sqrt{\frac{\sqrt{5}-1}{2}}<r<1$.
	\item[\rm 2.] $\frac{f}{g}\sim\frac{\sqrt{1+2\alpha}+1}{2\sqrt{\alpha}}(\sqrt{5}-1)^{\frac{d}{2}}\rightarrow\infty,$ if $r=\sqrt{\frac{\sqrt{5}-1}{2}}$.
	\item[\rm 3.] $\frac{f}{g}\sim\frac{1}{\sqrt{2}}(2\sqrt{1-r^2})^{d/2}\rightarrow\infty,$ if $0<r<\sqrt{\frac{\sqrt{5}-1}{2}}$.
\end{enumerate}
\end{Sl}

\textsc{Proof.} 
Obviously $f\sim \sqrt{\alpha 2^{d}}$ for $0\leq r <1.$

If $\sqrt{\frac{\sqrt{5}-1}{2}}<r<1,$ then $\frac{f}{g}\sim \frac{\sqrt{\alpha 2^{d}}}{\alpha/r^d}=
\frac{1}{\sqrt{\alpha}}(r\sqrt{2})^d\rightarrow\infty$ for $d\rightarrow\infty,$ since $r\sqrt{2}>1$ for $r>\sqrt{\frac{\sqrt{5}-1}{2}}.$

If $r=\sqrt{\frac{\sqrt{5}-1}{2}},$ then $\frac{f}{g}\sim  \frac{r^d(\sqrt{1+2\alpha}+1)\sqrt{\alpha 2^{d}}}{2\alpha}=\frac{\sqrt{1+2\alpha}+1}{2\sqrt{\alpha}}(r\sqrt{2})^d=\frac{\sqrt{1+2\alpha}+1}{2\sqrt{\alpha}}(\sqrt{5}-1)^{\frac{d}{2}}\rightarrow\infty$ for $d\rightarrow\infty.$ 

If $0<r<\sqrt{\frac{\sqrt{5}-1}{2}},$ then $\frac{f}{g}\sim \frac{\sqrt{\alpha 2^{d}}(1-r^2)^{d/4}}{\sqrt{2\alpha}}=
\frac{1}{\sqrt{2}}(2\sqrt{1-r^2})^{d/2}\rightarrow\infty$ for $d\rightarrow\infty,$ since $2\sqrt{1-r^2}>1$ for $0<r<\sqrt{\frac{\sqrt{5}-1}{2}}.$ The corollary is proved.

\section*{Acknowledgements}

Author is grateful to A.\,N.\,Gorban and N.\,Yu.\,Zolotykh for useful discussions. 
The work is supported by the Ministry of Education and Science of
Russian Federation (project 14.Y26.31.0022).

\end{document}